\documentclass[11pt]{article}

\usepackage{graphicx}
\usepackage{latexsym}
\usepackage[usenames,dvipsnames]{color}
\usepackage{amssymb,amsmath}
\usepackage{lscape}
\usepackage{verbatim}
\usepackage{bm}
\usepackage{enumitem} \setlist{nosep}
\usepackage{inputenc}

\normalbaselines
\catcode `\@=11
\@addtoreset{equation}{section}

\catcode`\@=12

\newtheorem{theorem}{Theorem}[section]
\newtheorem{defn}[theorem]{Definition}
\newtheorem{prop}[theorem]{Proposition}
\newtheorem{lemma}[theorem]{Lemma}
\newtheorem{cor}[theorem]{Corollary}

\newtheorem{rems}[theorem]{Remarks}
\newtheorem{rem}[theorem]{Remark}
\newtheorem{example}[theorem]{Example}

 \newcommand{\bea}{\begin{eqnarray}}
\newcommand{\ena}{\end{eqnarray}}

\newcommand{\beano}{\begin{eqnarray*}}
\newcommand{\enano}{\end{eqnarray*}}

\newcommand{\bei}{\begin{itemize}}
\newcommand{\eni}{\end{itemize}}

 \newcommand{\bee}{\begin{enumerate}}
 \newcommand{\ene}{\end{enumerate}}

\newcommand{\be}{\begin{equation}}
\newcommand{\en}{\end{equation}}

\newcommand{\bedefi}{\begin{defn} \rm }
\newcommand{\findefi}{\end{defn}}

\newcommand{\eofproof}{\hfill $\square$\bigskip}

\newcommand{\belem}{\begin{lemma}}
\newcommand{\enlem}{\end{lemma}}

\newcommand{\beprop}{\begin{prop}}
\newcommand{\enprop}{\end{prop}}

\newcommand{\betheo}{\begin{theorem}}
\newcommand{\entheo}{\end{theorem}}

\newcommand{\becor}{\begin{cor}}
\newcommand{\encor}{\end{cor}}

\newcommand{\berem}{\begin{rem} \rm}
\newcommand{\enrem}{\end{rem}}

\newcommand{\berems}{\begin{rems} \rm}
\newcommand{\enrems}{\end{rems}}

\newcommand{\beex}{\begin{example}$\!\!\!$ \rm }
\newcommand{\enex}{ \end{example}}

\newcommand{\norm}[2]{
\left\| #2 \right\|_{#1}
}
\newcommand{\Hil}[0]{\mathcal{H} }

 \newcommand{\ov}{\overline}

\newcommand{\nN}{\mathbb{N}}
\newcommand{\RN}{\mathbb{R}}

\newcommand{\ZN}{\mathbb{Z}}
\newcommand{\CN}{{\mathbb C}}

\newcommand{\bdim}{{\bf Proof. }}
 \newcommand{\edim}{\eofproof}

\newcommand{\hs}{Hilbert space}

\def\B{{\mathcal B}}

\def\D{{\mathcal D}}

\def\G{{\mathcal G}}

\def\H{{\mathcal H}}
\def\K{{\mathcal K}}

\def\P{{\mathcal P}}
\def\T{{\mathcal T}}

\newcommand{\NN}[0]{\mathbb{N}}
\newcommand{\ud}{\,\mathrm{d}}

\newcommand{\ip}[2]{ \langle {#1} |{#2}  \rangle}

\newcommand{\mc}{\mathcal}
\def\<{\langle}
\def\>{\rangle}

\definecolor{teal}{rgb}{0.0, 0.5, 0.5}

\makeatletter
\newcommand{\oset}[2]{%
	{\mathop{#2}\limits^{\vbox to -1\ex@{\kern-\tw@\ex@
				\hbox{\footnotesize  #1}\vss}}}}
\makeatother

\begin{document}

\begin{flushleft}
 {\LARGE   Lower semi-frames and metric operators} \vspace*{7mm}

{\large\sf   J-P. Antoine $\!^{\rm a}$, R. Corso $\!^{\rm b}$ and C. Trapani $\!^{\rm c}$
}
\\[3mm]
$^{\rm a}$  {\small Institut de Recherche en Math\'ematique et  Physique, Universit\'e catholique de Louvain \\
\hspace*{3mm}B-1348   Louvain-la-Neuve, Belgium\\
\hspace*{3mm}{\it E-mail address}: jean-pierre.antoine@uclouvain.be}
\\[1mm]
$^{\rm b,c}$ {\small Dipartimento di Matematica e Informatica,
Universit\`a degli Studi di Palermo, \\
\hspace*{3mm} I-90123 Palermo, Italy\\
$^{\rm b}$\hspace*{1.5mm}{\it E-mail address}: rosario.corso02@unipa.it}\\
$^{\rm c}${\small \hspace*{1.5mm}{\it E-mail address}: camillo.trapani@unipa.it}
\end{flushleft}

\textbf{Abstract}
 This paper deals with the possibility of transforming a weakly measurable function in a Hilbert space into a continuous frame by a metric operator, i.e., a strictly positive
self-adjoint operator. A necessary condition is that the domain of the analysis operator associated to the function be dense. The study is done also with the help of the generalized frame operator associated to a weakly measurable function, which has better properties than the usual frame operator. A special attention is given to lower semi-frames: indeed if the domain of the analysis operator is dense, then a lower semi-frame can be transformed into a Parseval frame with a (special) metric operator.

\section{Introduction}

In recent papers, one of us (RC)  \cite{corso,corso2} has analyzed sesquilinear forms defined by sequences in \hs s and operators associated to them by means of representation theorems.  In particular, he derived results about lower semi-frames and duality.

It turns out that most results from \cite{corso,corso2} can be extended to the continuous case and that is one of the aims of this present paper. The results are reported in Sect.3, but we give here a brief summary. The continuous case involves a locally compact space $(X,\mu)$ with a Radon measure $\mu$. A function $\phi:X\to \H, x\mapsto \phi_x$ is said to be {\it weakly measurable} if for every $f\in \H$ the function $x\mapsto \ip{f}{\phi_x}$ is measurable. A weakly measurable function $\phi$ is said to be \emph{$\mu$-total} if $\ip{f}{\phi_x}=0$ for a.e. $x\in X$ implies that $f=0$.
A weakly measurable function $\phi$ is a \emph{continuous frame} of $\H$ if there exist constants
{$0< {\sf m}\leq{\sf M}<\infty$ } (the  frame bounds) such that
\be\label{eq:frame}
{\sf m}  \norm{}{f}^2 \leq    \int_{X}  |\ip{f}{\phi_{x}}| ^2 \, \ud \mu(x)  \le {\sf M}  \norm{}{f}^2 , \qquad \forall \, f \in \H.
\end{equation}
If a weakly measurable function $\phi$ satisfies
\beano
\int_{X}  |\ip{f}{\phi_{x}}| ^2 \, \ud \mu(x)   \leq { \sf M}  \norm{}{f}^2 , \qquad \forall \, f \in \H, \,
\enano
then we say that $\phi$ is a {\it Bessel mapping}  of $\H$.
On the other hand, $ \phi$  is a \emph{lower semi-frame}  of $\H$ if  there exists   a constant ${\sf m}>0$ such that
\be
{\sf m}  \norm{}{f}^2 \leq  \int_{X}  |\ip{f}{\phi_{x}}| ^2 \, \ud \mu(x) , \qquad \forall \, f \in \H.
\label{eq:lowersf}
\end{equation}
Given a weakly measurable function $\phi$, the operator $C_\phi:D(C_\phi)\subseteq \H \to L^2(X,d\mu)$ with domain
$$D(C_\phi):=\left \{f\in \H: \int_{X}  |\ip{f}{\phi_{x}}| ^2 \, \ud \mu(x) <\infty \right\}$$
and $({C}_{\phi}f)(x) =\ip{f} {\phi_{x}}$, $f \in D(C_\phi),$
is called the {\it analysis} operator of $\phi$. We can define the sesquilinear form
$$
\Omega_\phi(f,g)=\int_X \ip{f}{\phi_x}\ip{\phi_x}{g}d\mu(x),\;\; f,g\in D(C_\phi)
$$
and associate a positive self-adjoint operator ${\sf T}_\phi$ in the space $\H_\phi$, { the closure of  $D(C_\phi)$ in $\H$}, by Kato's representation theorem \cite{kato}, which we call \emph{generalized frame operator}. When $\phi$ is a lower semi-frame of $\H$, then the range of ${\sf T}_\phi$  is $\H_\phi$ and the function $\psi:X\to \H$, defined by $\psi_x={\sf T}_\phi^{-1}P_\phi \phi_x, x\in X$, {where $P_\phi$ is the orthogonal projection onto $\H_\phi$,}
is a Bessel mapping and the reconstruction formula
$$
f=\int_X \ip{f}{\phi_{x}}\psi_x \ud\mu(x), \;\,
$$
holds for every $f\in D(C_\phi)$ in a weak sense.

In this paper we do not confine ourselves to extend results of \cite{corso,corso2}, but actually we set two more goals. From one hand, given a lower semi-frame $\phi:X\to \H$ with $D(C_\phi)$ dense in $\H$, we consider general powers ${\sf T}_\phi^{-k} \phi$ with $k\geq 0$. These functions are {Bessel mappings}, frames or lower semi-frames in the space $\H({\sf T}_\phi^m)$ (given by the domain of ${\sf T}_\phi^m$ and the inner product $\ip{{\sf T}_\phi^m \cdot}{{\sf T}_\phi^m \cdot}$) with $m\geq 0$ according to a simple relation between $k$ and $m$ (see Theorem \ref{th_T^k}).

When $\phi:X\to \H$ is a $\mu$-total weakly measurable function with $D(C_\phi)$ dense, then ${\sf T}_\phi$ is in particular a metric operator, i.e.,  a strictly positive self-adjoint operator. Metric operators are a topic familiar in the theory of the so-called  $\P\T$-symmetric quantum mechanics \cite{bender, mosta2}.  In our previous works
\cite{AT-metric1,AT-metric2,AT-metric3}, we have analyzed thoroughly the structure generated by such a metric operator, bounded or unbounded, namely a   lattice of \hs s.

As particular case of Theorem \ref{th_T^k}, if $\phi:X\to \H$ is a lower semi-frame with $D(C_\phi)$ dense, then ${\sf T}_\phi^{-1/2}\phi$ is a Parseval frame of $\H$. This inspired us to consider the following more general problem.

\smallskip
{\bf Question:} for which weakly measurable functions $\phi:X\to \H$ there exists a metric operator $G$ on $\H$ such that $\phi_x\in D(G)$ for all $x\in X$ and  $G\phi$ is a frame?
\smallskip

Partial answers to this problem are given in Theorem \ref{th_funct_metric}. In particular, necessary conditions are that $D(C_\phi)$ is dense, and that $\phi$ is $\mu$-total if $\phi$ is in addition a Bessel mapping. In the discrete case, if $\phi:\NN\to \H$ is a Schauder basis, then the problem has a positive solution (more precisely, one can again take $G=T^{-1/2}_\phi$ and $T^{-1/2}_\phi \phi$ is actually an orthonormal basis).

The paper is organized as follows.
After reviewing the conventional definitions about  frames and semi-frames in Sec.2, we introduce in Sec.3 the generalized frame operator ${\sf T}_\phi$, whose properties are more convenient that those of the standard frame operator $S_\phi$.
In Sec.4, we investigate the various (semi)-frames generated by a lower semi-frame. In Sec.5, we review the  lattice of \hs s generated by a metric operator. In Sec.6, we face the question of transforming functions in frames. We conclude in Sec.7 by several examples.

\section{Preliminaries}\label{sec-prel}

Before proceeding, we list further definitions and  conventions. The framework is
 a (separable) Hilbert space $\H$, with the inner product $\ip{\cdot}{\cdot}$ linear in the first factor.
$GL(\H)$ denotes the set of all invertible bounded operators on $\H$ with bounded inverse.
Throughout the paper, we will consider weakly measurable functions $\phi: X \to \H$, where $(X,\mu)$ is a  locally compact  space with a Radon measure $\mu$.

Given a continuous frame $\phi$, the  {analysis} operator is defined and bounded on $\H$, i.e. {${C}_{\phi}: \Hil \to L^{2}(X, \ud\mu)$}
\footnote{As usual, we identify a function $\xi$ with its residue class in  $L^{2}(X, \ud\mu)$.
}
and the corresponding \emph{synthesis operator} ${C}_{\phi}^\ast: L^{2}(X, \ud\mu) \to \H$ is defined as
 (the integral being understood in the weak sense,  as usual)
\be\label{eq:synthmap}
{C}_{\phi}^\ast \xi =  \int_X  \xi(x) \,\phi_{x} \; \ud\mu(x), \mbox{ for} \;\;\xi\in L^{2}(X, \ud\mu).
\end{equation}
Moreover, we set   $S_\phi:={C}_{\phi}^* {C}_{\phi}$, {which is self-adjoint.} Then it follows that
$$
\ip{S_\phi f}{g}=\ip{{C}^*_\phi {C}_\phi f}{g} = \ip{ {C}_\phi f}{{C}_\phi g} = \int_{X}  \ip{f}{\phi_{x}} \ip{\phi_x}{g} \, \ud \mu(x) .
$$
Thus, for continuous frames, $S_\phi$ and $S_\phi^{-1}$ are both bounded, that is, $S_\phi\in GL(\H)$.

Following \cite{ant-bal-semiframes1,ant-bal-semiframes2}, we will say that a function  $\phi$ is   a \emph{semi-frame} if
it   satisfies only one of the frame inequalities in \eqref{eq:frame}. We already introduced the lower semi-frames (if $\phi$ satisfies \eqref{eq:lowersf}, then it could still be a frame, hence we say that the lower semi-frame $\phi$ is \emph{proper} if it is not a frame). Note that the lower frame inequality automatically implies  that $\phi$ is  {$\mu$-total.}
On the other hand, a weakly measurable function $\phi$ is   an \emph{upper semi-frame} if is $\mu$-total, that is, $N(C_{\phi} )= \{0\}$,
and there exists ${\sf M}<\infty$ such that

\be\label{eq:upframe}
0 < \int_{X}  |\ip{f}{\phi_{x}}| ^2 \, \ud \mu(x)   \leq { \sf M}  \norm{}{f}^2 , \; \forall \, f \in \H, \, f\neq 0.
\en
{Thus an  upper   semi-frame is a total Bessel mapping \cite{gab-han}.}
Notice this definition does not forbid $\phi$ to be a frame. Thus we say that $\phi$ is a {\emph{proper} upper semi-frame} if it is not a frame.

If $\phi$ is a { proper upper semi-frame} $S_\phi$   is bounded and $S_\phi^{-1}$ is unbounded,
as follows immediately from \eqref{eq:upframe}.
{In the lower case, however, the definition of $S_\phi$  must be changed,
since the domain $D(C_\phi)$
need not be dense,    so that ${C}_{\phi}^*$ may not exist. Instead, following \cite[Sec.2]{ant-bal-semiframes1} one defines the synthesis operator  as
\be
{D}_{\phi}F =  \int_X  F(x) \,\phi_{x} \; \ud\mu(x),  \;\;F\in L^{2}(X, \ud\mu) ,
\label{eq:synthmap2}
\end{equation}
on the domain of all elements $F$ for which the integral in \eqref{eq:synthmap2} converges weakly in $\H$,
and  then {$S_\phi:={D}_{\phi} {C}_{\phi}$.}
With this definition, it is shown in \cite[Sec.2]{ant-bal-semiframes1} that if $\phi$ is a proper lower semi-frame then $S_\phi$  is unbounded and $S_\phi^{-1}$ is bounded.}
All these objects are studied in detail in our previous papers \cite{ant-bal-semiframes1,ant-bal-semiframes2}. In particular, it is shown
there that a natural notion of duality exists, namely,  two measurable functions $\phi, \psi$ are dual to each other (the relation is symmetric) if one has
\be\label{eq-dual}
\ip{f}{g} = \int_X  \ip{f}{\phi_{x}} \ip{\psi_{x}} {g}\, \; \ud\mu(x), \; \forall\, f,g \in \H.
\en
 This duality property   extends to  lower semi-frames and Bessel mappings, as shown in Proposition \ref{prop32} below.
 \smallskip

Consider the following sesquilinear form on the domain $\D_1 \times \D_2$ :
\be\label{eq:form}
\Omega_{\psi, \phi}(f,g) = \int_X \ip{f}{\psi_x} \ip{\phi_x}{g} \ud\mu(x), \; f\in \D_1, g\in \D_2.
\en
If $\D_1 =\D_2= \H$ and the form $\Omega_{\psi, \phi}$ is bounded on on $\H \times \H$, {that is, $| \Omega_{\psi, \phi}(f,g) | \leq c \norm{}{f}\norm{}{g}$,}   for some $c>0$, then
the couple of weakly measurable functions $(\psi, \phi)$   is called a \emph{reproducing pair} if
the corresponding bounded operator $S_{\psi, \phi}$ given weakly  by
$$
\ip{S_{\psi, \phi}f}{g} = \int_X \ip{f}{\psi_x} \ip{\phi_x}{g}  \ud\mu(x) , \; \forall\,f, g\in\H,
$$
belongs to $GL(\H)$. If $\psi = \phi$, we recover the notion  of continuous frame.

\medskip

Under certain conditions, boundedness of $\Omega_{\psi, \phi}$ is automatic, as shown in \cite[Prop. 7]{corso}
\beprop
\label{prop21}
If $\D_1 =\D_2= \H$, $X$ is  locally compact and  $\sigma$-compact, that is, $X = \bigcup_{n} K_{n}, K_{n} \subset  K_{n+1}$, with $K_{j}$
compact for every $j$,
and
$\sup_{x\in X}(\norm{\H}{\phi_x }\norm{\H}{\psi_x }) <\infty$,
then the form $\Omega_{\psi, \phi}$ is bounded  on $\H \times \H$.
\enprop
\bdim
Define
$$
\Omega^n_{\psi, \phi}(f,g) := \int_{K_n}  \ip{f}{\psi_x} \ip{\phi_x}{g} \ud\mu(x), f,g \in\H.
$$
By assumption, there exists $c>0$ such that
$$
|\Omega^n_{\psi, \phi}(f,g)| \leq c \norm{}{f} \norm{}{g} \left | \int_{K_n} \ud\mu(x) \right| < \infty
$$
Hence there exists a bounded operator $T_n$ such that $\Omega^n_{\psi, \phi}(f,g) = \ip{T_n f}{g}$.
Applying the Banach-Steinhaus theorem to the functional $g\mapsto \Omega_{\psi, \phi}(f,g) $, one gets that the operator $T_{\psi, \phi}$ associated to
$\Omega_{\psi, \phi}$ is defined on the whole of $\H$. Doing the same with $T_n^\ast$, as in  \cite[Prop. 7]{corso}, we conclude that the form
$\Omega_{\psi, \phi}$ is bounded  on $\H \times \H$.
\edim

	The converse of Proposition \ref{prop21} does not hold, as shown in the following example (which is based on \cite[Example 2.5]{balaszetal}), {in the sense that boundedness of the form   $\Omega$ does not imply the supremum condition.}
	
	\beex
	Let $h\in \H\backslash\{0\}$ and let $a:\RN\to \CN$ be such that $a\in L^2(\RN)\backslash L^\infty(\RN)$.  Define $\phi_x=a(x)h$ for all $x\in \RN$. Then $\phi:\RN \to \H$ is a weakly-measurable function and a Bessel mapping because
	$$
	\int_\RN |\ip{f}{\phi_x}|^2 dx \leq \norm{L^2(\RN)}{a}^2 \|h\|^2\|f\|^2, \quad  f\in \H.
	$$
	In conclusion, the sesquilinear form $\Omega_{\psi, \phi}$, where $\psi=\phi$, is defined and bounded on $\H\times \H$, but  $\sup_{x\in X}(\norm{\H}{\phi_x }\norm{\H}{\psi_x }) =\infty$.
	\enex

\section{The generalized frame operator ${\sf T}_\phi$}
\label{sec-Tphi}

In the previous section we defined the {\it frame operator} $S_\phi$ for a lower semi-frame. However, this operator lacks  good properties, in general (for instance $S_\phi$ need not be self-adjoint like in the case of an upper semi-frame, even if $S_\phi$ is non-negative). In this section, we are going to construct a new operator associated to $\phi$ which plays the r\^{o}le of $S_\phi$ for lower semi-frames. We show its main properties, in particular, concerning the definition of a Bessel dual mapping in a natural way.

We note that if $\phi$ is a  proper lower semi-frame, the r.h.s. of \eqref{eq:lowersf}  actually diverges for some $f$.
As already said, the domain $D(C_\phi)$ need not be dense in $\H$. It is useful to work with the Hilbert space $\H_\phi$ made of the closure of $D(C_\phi)$ endowed with the topology of $\H$.
	
	The analysis operator $C_\phi$ is closed \cite[Lemma 2.1]{ant-bal-semiframes1}. Therefore the sesquilinear form
	$$
	\Omega_\phi(f,g)={ \ip{ {C}_\phi f}{{C}_\phi g}}=\int_X \ip{f}{\phi_x}\ip{\phi_x}{g}d\mu(x),\;\; f,g\in D(C_\phi)
	$$
	is non-negative and closed. By Kato's second representation theorem \cite[Theorem 2.23]{kato} there exists an operator ${\sf T}_\phi:D({\sf T}_\phi)\subset \H_\phi\to \H_\phi $, {with $D({\sf T}_\phi)\subset D(C_\phi)$}, such that\footnote{We use this sans serif font in order to avoid confusion with the generalized synthesis operator
	$T_\phi$ introduced in our papers about reproducing pairs \cite{jpa_ms_ct,at-reprodpairs2,AT-axioms}.}
	\begin{itemize}
		\item $\Omega_\phi(f,g)=\ip{{\sf T}_\phi f}{g}$ for all $f\in D({\sf T}_\phi)$ and $g\in D(C_\phi)$;
		\item if $T:D(T)\subset \H_\phi\to \H_\phi$  is such that $\Omega_\phi(f,g)=\ip{T f}{g}$ for all $f\in D(T)$ and $g\in D(C_\phi)$, then $T\subset {\sf T}_\phi$;
		\item ${\sf T}_\phi$ is non-negative and self-adjoint in $\H_\phi$;
		\item $\Omega_\phi(f,g)=\ip{{\sf T}_\phi^{1/ 2} f}{{\sf T}_\phi^{1/2} g}$, for all $f,g\in D(C_\phi)$.
	\end{itemize}
	
	We call ${\sf T}_\phi$ the {\it generalized frame operator} of $\phi$. The motivation behind this name is that when $\phi$ is a continuous frame then ${\sf T}_\phi=S_\phi$. The generalized frame operator has been studied in \cite{corso,corso2} in the discrete case and a preliminary extension to the continuous setting has been given in	\cite{bellomonte}.
	
	If $D(C_\phi)$ is dense, then of course $\H_\phi=\H$ (i.e., ${\sf T}_\phi$ is a non-negative self-adjoint operator on $\H$) and $S_\phi \subset {\sf T}_\phi$.
	In \cite{corso2}, it was proved that in the discrete setting ($X=\NN$) we may have a strict inclusion $S_\phi\subsetneq {\sf T}_\phi $ (recall that in our paper $S_\phi$ is weakly defined, therefore in the discrete setting it corresponds to the operator $W_\phi$ in \cite{corso2}).
	
	Since $\Omega_\phi(f,g)=\ip{C_\phi f}{C_\phi g}$ we have ${\sf T}_\phi:=|C_\phi|^2=C_\phi^\times C_\phi$ where $C_\phi^\times$ and $|C_\phi|$ are the adjoint and the modulus of $C_\phi$ when we think of it as an operator $C_\phi:\H_\phi\to L^2(X,\mu)$.

	 The following characterization can be proved as in \cite{corso}.

	\beprop \label{prop31}
	Let $\phi$ be a weakly measurable function and ${\sf m}>0$. The following statements are equivalent.
	
	(i) $\phi$ is a lower semi-frame of $\H$ with lower bound {\sf m};
	
	(ii)  $\Omega_\phi$ is bounded from below by {\sf m}, i.e.,
	$$
	\Omega_\phi(f,f)\geq {\sf m}\|f\|^2,\;\;  \forall f\in D(C_\phi);
	$$

	(iii)  $C_\phi$ is bounded from below by $\sqrt{\sf m}$, i.e.
	$$
	\|C_\phi f\|\geq \sqrt {\sf m}\|f\|,\;\;  \forall f\in D(C_\phi);
	$$
	
	(iv)  ${\sf T}_\phi$ is bounded from below by {\sf m}, i.e.,
	$$
	\|{\sf T}_\phi f\|\geq {\sf m}\|f\|,\;\;  \forall f\in D({\sf T}_\phi);
	$$
	
	(v)  ${\sf T}_\phi$ is invertible and ${\sf T}_\phi^{-1}\in \B(\H_\phi)$ with $\|{\sf T}_\phi^{-1}\|\leq {\sf m}$.
	
	\enprop
	\bdim
		(i) $\iff$ (ii) {and (i) $\iff$ (iii)} By definition.
	
		(ii) $\implies$ (iv)
	Let $f\in D({\sf T}_\phi)$. Then $\|{\sf T}_\phi f\|\|f\|\geq \ip{{\sf T}_\phi f}{f}\geq {\sf m}\|f\|^2$. Thus $\|{\sf T}_\phi f\|\geq {\sf m}\|f\|$.

	(iv) $\implies$ (ii)  For $f\in D({\sf T}_\phi)$ we have $\Omega_\phi(f,f)=\ip{{\sf T}_\phi^{1/2} f}{{\sf T}_\phi^{1/2} f}\geq {\sf m}\|f\|^2$ by hypothesis. The inequality now extends to every $f\in D(C_\phi)$ noting that $D({\sf T}_\phi)$ is a core of $\Omega_\phi$ (see \cite[Theorem 2.1]{kato}).

	(iv) $\iff$ (v) This is clear since ${\sf T}_\phi$ is self-adjoint.
	
	\edim
	
	Now assume that $\phi$ is a lower semi-frame of $\H$. By \cite[Proposition 2.6]{ant-bal-semiframes1} there exists a Bessel mapping
	$\psi:X\to \H$ dual to $\phi$, i.e. for some ${\sf M}>0$
	$$
	 \int_{X}  |\ip{f}{\psi_{x}}| ^2 \, \ud \mu(x)  \le {\sf M}  \norm{}{f}^2 ,  \;\;\forall \, f \in \H,
	$$
	and
	\begin{equation}
	\label{dual_lower}
	f=\int_X \ip{f}{\phi_{x}}\psi_x \ud\mu(x), \;\, \text{for all } f\in D(C_\phi) \text{ in weak sense. }
	\end{equation}
 A proof may be found in \cite[Lemma 2.5 and Proposition 2.6]{ant-bal-semiframes1}. 	Note that the proof given there is incomplete, in the sense there is no guarantee that $\psi$ is total, see Item (4) in Proposition \ref{prop32} below.
Here we give a different proof involving the operator ${\sf T}_\phi$ using the same argument as in \cite[Theorem 4.1]{corso2}.

	\beprop
	\label{prop32}
	Let $\phi$ be a lower semi-frame of $\H$.
	Let $\psi$ be the function defined by  $\psi={\sf T}_\phi^{-1}P_\phi \phi$, i.e.  $\psi_x={\sf T}_\phi^{-1}P_\phi \phi_x$ for $x\in X$, where $P_\phi:\H\to \H_\phi$ is the orthogonal projection of $\H$ onto $ \H_\phi$. The following statements hold.
	\begin{itemize}
		\item [(1)] $\psi$ is a weakly measurable function;
		\item[(2)]  $\psi$ a Bessel mapping of $\H$;
		\item[(3)]  $\psi$ is dual to $\phi$, i.e. satisfies \eqref{dual_lower}.
		\item[(4)]  If $D(C_\phi)$ is dense, then $\psi$ is an upper semi-frame of $\H$.
	\end{itemize}
	\enprop
	\bdim
	We recall that ${\sf T}_\phi$ is invertible with ${\sf T}_\phi^{-1}\in \B(\H_\phi)$.
		
		(1) This is immediate.
		
		(2)
		Let $f\in \H$. Then
		\bea
		\nonumber \int_{X} |\ip{f}{{\sf T}_\phi^{-1}P_\phi \phi_x}|^2 \ud\mu(x)
		&=&\int_{X} |\ip{P_\phi f}{{\sf T}_\phi^{-1}P_\phi \phi_x}|^2\ud\mu(x)=\int_{X} |\ip{{\sf T}_\phi^{-1} P_\phi f}{ \phi_x}|^2\ud\mu(x)\\
		&=& \| {\sf T}_\phi^{1/2} {\sf T}_\phi^{-1} P_\phi f \|^2= \|{\sf T}_\phi^{-1/2} P_\phi f \|^2  {\leq}
		\|{\sf T}_\phi^{-1/2}\|^2 \|f \|^2. \label{eq1}
		\ena
		  		
		  (3) Let $f\in D(C_\phi)$ and $h\in \H$. Then
		  \beano
		  \ip{f}{h}&=&\ip{f}{P_\phi h}=\ip{f}{{\sf T}_\phi {\sf T}_\phi^{-1}P_\phi h}=\int_X  \ip{f}{\phi_x}\ip{\phi_x}{{\sf T}_\phi^{-1}P_\phi h}\ud\mu(x)\\
		  &=&\int_X  \ip{f}{\phi_x}\ip{P_\phi \phi_x}{{\sf T}_\phi^{-1}P_\phi h}\ud\mu(x)= \int_X  \ip{f}{\phi_x}\ip{{\sf T}_\phi^{-1}P_\phi \phi_x}{h}\ud\mu(x).
		  \enano
		  		
		  (4) This is immediate from \eqref{dual_lower}.
	\edim

{\berem
Part (3) of the above proposition comes from general properties of semi-bounded operators. If $\H, {\mc K}$ are \hs s,  and $C:D(C)\subset \H \to {\mc K}$ is a closed  operator such that $D(C)$ is dense and $\|Cf\|\geq \gamma \|f\|$, for every $f\in D(C)$, then by Kato's theorem there exists $T\geq 0$, with the properties described at the beginning of Section \ref{sec-Tphi}, such that $\ip{Cf}{Cg}= \ip{T^{1/2}f}{T^{1/2}g}$ for every $f,g \in D(C)$; as before $T$ is invertible with bounded inverse. On the other hand,
$$
|\ip{f}{g}| \leq \frac{1}{\gamma^2} \|f\|_C \|g\|_C, \quad \forall f,g \in D(C), \mbox{where} \;\|f\|_C: = \norm{}{Cf}.
$$
Hence, there exists an operator $X$, bounded in $\H(C):=D(C)[\|\cdot\|_C]$, such that
$$
 \ip{f}{g}= \ip{Cf}{CXg}, \quad \forall f,g \in D(C).
 $$
 Then,
$$
\ip{Cf}{CXg}= \ip{T^{1/2}f}{T^{1/2}Xg}=\ip{f}{g}, \quad \forall f,g \in D(C).
$$
This implies that $T^{1/2}Xg \in D(T^{1/2})=D(C)$, for every $g\in D(C)$ and $TXg=g$, $\forall g\in D(C)$. Thus $X\subset T^{-1}$.
In conclusion, $\ip{f}{g}=\ip{Cf}{CT^{-1}g}$, for all $f,g \in D(C)$. So, if $C=C_\phi$ is densely defined, the {\em dual} is found. If $D(C_\phi)$ is not dense we can proceed with the projection $P_\phi$ as before.
\enrem
}
Finally, as mentioned in \cite{corso2} for the discrete case,  calculations similar to \eqref{eq1} show the next result.
	
	\begin{prop}
		\label{prop 34}
		Let $\phi$ be a lower semi-frame of $\H$. Then the function $\psi= {\sf T}_\phi^{-1/2}P_\phi \phi$ is a Parseval frame for $\H_\phi$.		
	\end{prop}

Thus, following the standard terminology in frame theory, we can call ${\sf T}_\phi^{-1}P_\phi \phi$ and ${\sf T}_\phi^{-1/2}P_\phi \phi$ the {\it canonical dual Bessel mapping} and the {\it canonical tight frame} of $\phi$.

Now we exploit Proposition 5 of \cite{corso} and the discussion after it.
\beprop \label{prop5}
Let $\phi$ be a weakly measurable function of $\H$. Then  $\phi$ is a lower semi-frame of $\H$ if and only if there exists an inner product $\ip{\cdot}{\cdot}_+$
inducing a norm $\norm{+}{\cdot}$ on $D(C_\phi)$ for which $ D(C_\phi)[\norm{+}{\cdot}]$  is complete, continuously embedded into $\H$ and for some {\sf m}, {\sf M} $> 0$ one has
$$
{\sf m}  \norm{+}{f}^2  \leq   \int_{X}  |\ip{f}{\phi_{x}}| ^2 \, \ud \mu(x)  \leq {\sf M}  \norm{+}{f}^2 ,  \; \forall \, f \in D(C_\phi).
$$
\enprop
\bdim
It is sufficient to take
$
\norm{+}{f}^2 = \norm{C_\phi}{f}^2= \int_{X}  |\ip{f}{\phi_{x}}| ^2 \, \ud \mu(x),
$
 for $f\in  D(C_\phi)$.
\edim

Let $\phi$ be a lower semi-frame in the \hs\ $\H$, with domain $D(C_\phi)$, assumed to be dense.
For every $x\in X$, the map $f\mapsto \ip{f}{\phi_x}$ is a bounded linear functional on the \hs\    {$ \H(C_\phi):= D(C_\phi)[\norm{+}{\cdot}]$.} By the Riesz Lemma, there exists  an element $\chi_x^{\phi} \in D(C_\phi)$ such that
\be \label{eq8}
\ip{f}{\phi_x} = \ip{f}{\chi_x^{\phi}}_+  \,\,\, \forall  \, f \in D(C_\phi).
\en
 By Prop. \ref{prop5}, $\chi^{\phi}$ is a frame.

   We can explicitly determine the element $\chi^{\phi}$ when $\norm{+}{\cdot}$ is  the norm $\norm{1/2}{f}^2 =
    \norm{}{{\sf T}_\phi^{1/2}f}^2 $.
  Notice that, by {Prop. \ref{prop31}$(iv),$ }  this   norm   is equivalent to the the graph norm
 of ${\sf T}_\phi^{1/2}$.

  Then we have $\ip{f}{\phi_x} =\ip{f}{\chi_x^{\phi}}_{1/2}=\ip{ f}{{\sf T}_\phi\chi_x^{\phi}}$ for all $f\in D(C_\phi)$. Thus $\chi_x^{\phi}={\sf T}_\phi^{-1} \phi_x$ for all $x\in X$, i.e. $\chi^\phi$  is the canonical dual Bessel mapping of $\phi$.

 Following
the notation of  \cite{AT-PIPmetric},  denote by $\H( {\sf T}_\phi^{1/2})$ the \hs\  $ D({\sf T}_\phi^{1/2})$  with the   norm  $\norm{1/2}{\cdot}$. {  In the same way, denote by  $\H( C_\phi)$ the \hs\  $ D(C_\phi)$  with inner product  $\ip{C_\phi \cdot}{C_\phi \cdot}$}.  Hence we have proved the following result.

  \beprop
  	Let $\phi$ be a lower semi-frame of $\H$ with $D(C_\phi)$ dense. Then the canonical dual Bessel mapping of $\phi$ is a tight frame for the Hilbert space
  	$\H(C_\phi)= \H( {\sf T}_\phi^{1/2})$.
  \enprop

We can proceed in the converse direction, i.e. starting with a frame $\chi \in D(C_\phi)$, does there exist a lower semi-frame $\eta$  of $\H$ such that $\chi $ is the frame $\chi^{\eta}$ constructed from $\eta$ in the  way described above? The answer is formulated in the following

\beprop[{\cite[Prop. 6]{corso}}]
Let $\chi $ be a frame of $\H(C_\phi)=\H({\sf T}_\phi^{1/2} )$. Then
\\
(i) There exists a lower semi-frame $\eta$ of $\H $ such that $\chi  = \chi^{\eta}$  if, and only if, $\chi  \in D({\sf T}_\phi)$.
\\
(ii) If $\chi  =\chi^{\eta}$ for some lower semi-frame $\eta$ of $\H $, then $\eta= {\sf T}_\phi\chi$.
\enprop

\section{The functions generated by a lower semi-frame}
\label{sec-powers_low}

Throughout this section we continue to consider   a lower semi-frame $\phi$ in $\H$, with  $D(C_\phi)$ dense.

Since ${\sf T}_\phi^{-1}$ is defined on $\H$ we can actually apply different powers of ${\sf T}_\phi^{-1}$ on $\phi$ and get the functions ${\sf T}_\phi^{-k} \phi$, $k\in [0,\infty)$. Hence we can ask  for the properties of ${\sf T}_\phi^{-k} \phi$. Of course the answer depends on the Hilbert space where ${\sf T}_\phi^{-k} \phi$ is considered. For instance, for $k=0$ we have a lower semi-frame of $\H$ and, as seen in the previous section, for $k=\frac{1}{2}$ we have a frame for $\H$, while for $k=1$ we have a frame for $\H( {\sf T}_\phi^{1/2})$.
When we have powers of an unbounded,  closed, densely defined operator, then we can consider {\it scales} and {\it lattices} of Hilbert spaces, which we will consider in more detail in Section \ref{sec-metric}.
{  For a while, let us simply denote  by $\H({\sf T}_\phi^m)$, $m\geq 0$, the domain of ${\sf T}_\phi^m$ considered as a Hilbert space with norm $\|f\|_m= \|{\sf T}_\phi^m f\|$, $f\in  D({\sf T}_\phi^m)$. Then, if $m>n\geq 0$, we have $\H({\sf T}_\phi^m) \subset \H({\sf T}_\phi^n) \subset \H$.
The function ${\sf T}_\phi^{-k} \phi$ is not always a well-defined function of  $\H({\sf T}_\phi ^{m})$;
a sufficient condition for  ${\sf T}_\phi^{-k} \phi\in \H({\sf T}_\phi ^{m}$) is that $k\geq {m}$.

Having at our disposal the notion of scale of Hilbert spaces, we now come back to a lower semi-frame $\phi$ with $D(C_\phi)$ dense and to the functions ${\sf T}_\phi^{-k}\phi$ with $k\geq 0$.  For simplicity of notation, we write in a compact way the inner product of $\H({\sf T}_\phi^m)$ in the following way
 $\ip{f}{g}_m:=\ip{f}{g}_{{\sf T}_\phi^m}=\ip{{\sf T}_\phi^m f}{{\sf T}_\phi^m g}$.

\betheo
\label{th_T^k}
	Let $\phi$ be a lower semi-frame of $\H$ with $D(C_\phi)$ dense and $m,k\in [0,\infty)$, $k\geq m$. Then the following statements hold.
	\begin{enumerate}
		\item[(i)]  ${\sf T}_\phi^{-k}\phi $ is a Bessel mapping of $\H({\sf T}_\phi^m)$ if and only if $k\geq m+\frac{1}{2}$;
		\item[(ii)] ${\sf T}_\phi^{-k}\phi $ is a lower semi-frame of $\H({\sf T}_\phi^m)$ if and only if $m\leq k \leq m+\frac{1}{2}$;
		\item[(iii)] ${\sf T}_\phi^{-k}\phi $ is a frame of $\H({\sf T}_\phi^m)$ if and only if ${\sf T}_\phi^{-k}\phi $ is a Parseval frame of $\H({\sf T}_\phi^m)$, if and only if $k= m+\frac{1}{2}$.
	\end{enumerate}
\entheo
	\bdim
	We have for $f\in \H({\sf T}_\phi^m)$
	\bea\int_{X} \nonumber |\ip{f}{{\sf T}_\phi^{-k} \phi_x}|_{m}^2 \ud\mu(x)
	&=&\int_{X} |\ip{{\sf T}_\phi^{m} f}{{\sf T}_\phi^{m}{\sf T}_\phi^{-k} \phi_x}|^2\ud\mu(x)=\int_{X} |\ip{{\sf T}_\phi^{2m-k}  f}{ \phi_x}|^2\ud\mu(x)\\
	&=& \|{\sf T}_\phi^{1/2} {\sf T}_\phi^{2m-k} f \|^2= \| {\sf T}_\phi^{2m-k+1/2} f \|^2. \label{Pars_T^k}
	\ena
	Hence, taking into account that $0\in \rho({\sf T}_\phi)$ and that $	\|f\|^2_{{\sf T}_\phi^m}=\|{\sf T}_\phi^m f\|^2	$, (i) and (ii) follow by comparing $2m-k+\frac{1}{2}$ and $m$. We give the details for (i) as example.\
	
	Assume that  ${\sf T}_\phi^{-k}\phi $ is a Bessel mapping of $\H({\sf T}_\phi^m)$, then there exists $B>0$ such that for all $f\in \D({\sf T}_\phi^m)$ we have $\| {\sf T}_\phi^{2m-k+1/2} f \|^2\leq B \|{\sf T}_\phi^m f\|^2$. Rewriting the inequality as $\| {\sf T}_\phi^{m-k+\frac{1}{2}} {\sf T}_\phi^{m} f \|^2\leq B \|{\sf T}_\phi^m f\|^2$, we find that ${\sf T}_\phi^{m-k+\frac{1}{2}}$ is bounded,  i.e. $m-k+\frac{1}{2}\leq 0$. The other implication trivially holds by the same estimate.
	
	Finally, combining the two cases above, we obtain that  ${\sf T}_\phi^{-k}\phi $ is a frame of $\H({\sf T}_\phi^m)$ if and only if $k= m+\frac{1}{2}$. Moreover, in this case ${\sf T}_\phi^{-k}\phi $ is actually a Parseval frame, as one can see in \eqref{Pars_T^k}.
	\edim

\medskip

In particular, we recover some cases we discussed in the previous section, namely  $(k,m)=(1,\frac{1}{2})$ and $(k,m)=(1,0)$  (corresponding to the canonical Bessel mapping of $\phi$) and $(k,m)=(\frac{1}{2},0)$ (corresponding to the canonical tight frame of $\phi$).

{It is possible to generalize   Theorem \ref{th_T^k} by considering more general functions than powers of ${\sf T}_\phi$. For instance we could take   the set $\Sigma$ of real valued functions $\sf g$ defined on the spectrum $\sigma({\sf T}_\phi)$, which are measurable with respect to the spectral measure of ${\sf T}_\phi$ and such that ${\sf g}$ and $\widetilde {\sf g}: =1/{\sf g}$ are bounded on compact subsets of $\sigma({\sf T}_\phi)$. For every ${\sf g} \in \Sigma$ we denote by $\H_{\sf g}$ the Hilbert space completion of $D({\sf g}({\sf T}_\phi))$  with respect to the norm $\|f\|_{\sf g}= \|{\sf g}({\sf T}_\phi)f\|$, $f\in D({\sf g}({\sf T}_\phi))$. As shown in \cite[Sec.10.4]{ait_book},
we get a LHS if the order is defined by ${\sf h} \preceq {\sf g} \;\Longleftrightarrow \;\exists \;\gamma>0$ such that ${\sf h} \leq \gamma {\sf g}$. We put ${\sf i}(t):=t, \;   t\in \sigma({\sf T}_\phi)$.  Then we get
\betheo
\label{th_52}
	Let $\phi$ be a lower semi-frame of $\H$ with $\D(C_\phi)$ dense and ${\sf g},{\sf h}$ non-negative functions from $\Sigma$, with ${\sf g}\succeq {\sf h}$. Suppose that $\widetilde{\sf g}$ and ${\sf h}\,\widetilde{{\sf g}}$ are bounded functions.
	 Then the following statements hold.
	\begin{enumerate}
		\item[(i)]  $\widetilde{{\sf g}}({\sf T}_\phi)\phi $ is a Bessel mapping of $\H_{\sf h}$ if and only if ${\sf i}^{1/2}{\sf h}\preceq {\sf g}$;
		\item[(ii)] $\widetilde{{\sf g}}({\sf T}_\phi)\phi $ is a lower semi-frame of $\H_{\sf h}$ if and only if ${\sf h} \preceq {\sf g}\preceq {\sf i}^{1/2}{\sf h}$;
		\item[(iii)] $\widetilde{{\sf g}}({\sf T}_\phi)\phi $ is a frame of $\H_{\sf h}$ if and only if $\widetilde{{\sf g}}({\sf T}_\phi)\phi$ is a Parseval frame of $\H_{\sf h}$, if and only if ${\sf i}^{1/2}{\sf h}\preceq {\sf g}$.
	\end{enumerate}
\entheo
\bdim The proof is similar to that of Theorem \ref{th_T^k}. In fact, for $f\in \H_{\sf h}$, using the functional calculus for ${\sf T}_\phi$, we have
\bea\int_{X} \nonumber |\ip{f}{\widetilde{{\sf g}}({\sf T}_\phi) \phi_x}|_{\sf h}^2 \ud\mu(x)
	&=&\int_{X} |\ip{{\sf h}({\sf T}_\phi)f}{({\sf h}\widetilde{\sf g})({\sf T}_\phi) \phi_x}|^2\ud\mu(x)=\int_{X} |\ip{({\sf h}^2\widetilde{\sf g})({\sf T}_\phi) f}{ \phi_x}|^2\ud\mu(x)\\
	&=& \|{\sf T}_\phi^{1/2} ({\sf h}^2\widetilde{\sf g})({\sf T}_\phi) f \|^2= \| ({\sf h}^2\widetilde{\sf g}{\sf i}^{1/2})({\sf T}_\phi) f \|^2. \label{Pars_T^k2}
	\ena
Thus, for instance, if ${\sf i}^{1/2}{\sf h}\preceq {\sf g}$, then ${\sf i}^{1/2}{\sf h}^2\widetilde{\sf g}\preceq {\sf h}$; hence $\widetilde{{\sf g}}({\sf T}_\phi)\phi $ is a Bessel mapping of $\H_{\sf h}$.
The rest of the proof is analogous.
\edim}
}
\section{Metric operators}
\label{sec-metric}
The generalized frame operator of a total weak measurable function $\phi$ with $D(C_\phi)$ dense is an example of metric operator in the sense of the following definition
\cite{AT-metric1}.
{
\bedefi By a metric operator in a Hilbert space $\H$, we mean a strictly positive
self-adjoint operator $G$, that is, $G > 0$ or $\ip{Gf}{f} \geq 0$ for every $f \in D(G)$ and $\ip{Gf}{f} =
0$  if and only if $f= 0$.
\findefi
}

Let $S$ be an unbounded  closed operator with dense domain  $D(S)$. As usual, define the graph norm of  $S$:
\begin{align*}
\ip{f}{g}_{\rm gr} &:= \ip{f}{g} + \ip{Sf}{Sg}, \; f,g \in D(S),
\\
\norm{\rm gr}{f}^2 &=\norm{}{f}^2  + \norm{}{Sf}^2 .
\end{align*}
Then the norm  $\norm{\rm gr}{\cdot}$ makes $D(S)$ into a \hs\ continuously embedded into  $\H$.
For $f\in  D(S^\ast S)$, we may write
$\norm{{\rm gr}}{f}^2 =   \ip{f}{(I+ S^\ast S)f}$.
Note that the operator $ S^\ast S = |S|^2$ is self
adjoint and non-negative: $|S|^2 \geq 0$. In addition, $D(|S|) = D(S)$ and $ N(S^\ast S) = N(S)$ \cite[Theor.5.39 and 5.40]{weidmann}.

Given $S$ as above, the operator $R_S:= I+ S^\ast S$
is self-adjoint, with domain  $D(S^\ast S)$,
 and   $R_S\geq 1$.
  Hence $R_S$ is an unbounded metric operator, with bounded inverse
$R_S^{-1}=( I+ S^\ast S) ^{-1}$.
{In our previous works \cite{AT-metric1,AT-metric2,AT-metric3}, we have analyzed the lattice of \hs s generated by such a metric operator. In the sequel, we summarize this discussion, keeping the same notations.}

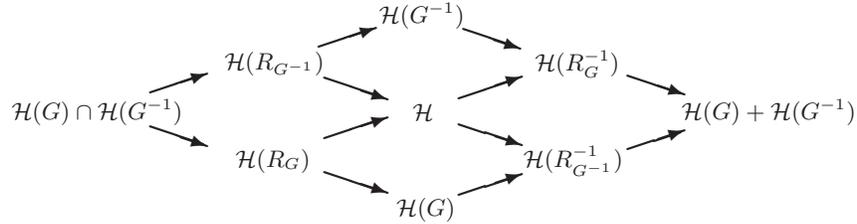
\begin{figure}[h]
\centering \setlength{\unitlength}{0.38cm}
\begin{picture}(8,8)

\put(3.5,4){
\begin{picture}(8,8) \thicklines
\footnotesize
 \put(-3.4,-0.9){\vector(3,1){2.2}}
\put(-3.6,2.3){\vector(3,1){2}}
 \put(-3.4,-2.1){\vector(3,-1){2.2}}
\put(-3.4,1.2){\vector(3,-1){2.2}}
\put(1.3,0.5){\vector(3,1){2.2}}
\put(1.3,-2.9){\vector(3,1){2.2}}
 \put(1.3,-0.4){\vector(3,-1){2.2}}
\put(1.45,2.9){\vector(3,-1){2}}
\put(0,3.4){\makebox(0,0){$ \H(G^{-1})$}}
\put(-0.1,0){\makebox(0,0){ $\H$}}
\put(0,-3.4){\makebox(0,0){ $\H(G)$}}
\put(-11.5,0){\makebox(0,0){ $\H(G) \cap  \H(G^{-1})$}}
\put(12,0){\makebox(0,0){ $\H(G) +  \H(G^{-1})$}}
\put(-9.5,0.7){\vector(3,1){2}}
\put(7.2,-1.3){\vector(3,1){2}}
\put(7.2,1.3){\vector(3,-1){2}}
\put(-9.5,-0.5){\vector(3,-1){2}}
\put(-5.3,1.7){\makebox(0,0){ $\H(R_{G^{-1}})$}}
\put(-5.2,-1.7){\makebox(0,0){$\H(R_{G})$}}
\put(5.3,1.7){\makebox(0,0){ $ \H(R_G^{-1})$}}
\put(5.3,-1.7){\makebox(0,0){$ \H(R_{G^{-1}}^{-1})$}}

\end{picture}
}
\end{picture}
\caption{\label{fig:diagram}The lattice of \hs s generated by a metric operator.}

\end{figure}

In the general case where both the metric operator $G$ and its inverse $G^{-1}$ are unbounded, the lattice is given in Fig.  \ref{fig:diagram}.
Given the metric operator $G$, equip the domain $D(G^{1/2})$
with the following norm
\be\label{norm-RG}
\norm{R_G}{f}^2 =  \norm{}{(I+G)^{1/2}f}^2,\; f \in D(G^{1/2}),
\en
 Since this norm is  equivalent to the graph norm of  $G^{1/2}$,
this makes $ D(G^{1/2})$ into a \hs, denoted $\H(R_G)$, dense in $\H$.
 Next,  we equip $\H(R_G)$ with the norm $\norm{G}{f}^2 := \norm{}{G^{1/2}f}^2$
and denote by $\H(G)$ the completion of $\H(R_G)$ in that norm and corresponding inner product $\ip{\cdot}{\cdot}_G :=\ip{G^{1/2}\cdot}{G^{1/2}\cdot}$.  Hence, we have  $\H(R_G) = \H \cap \H(G)$, with the so-called projective norm \cite[Sec.I.2.1]{pip-book}, which here is simply the graph norm of
$G^{1/2}$.

Next we proceed in the same way with the inverse operator $G^{-1}$, and we obtain another \hs, $\H(G^{-1})$. Then  we consider the lattice generated by
$\H(G)$ and  $ \H(G^{-1})$ with the operations
 \begin{align}\label{eq:lattice1}
 \H_1 \wedge \H_2& := \H_1 \cap \H_2 \, , \\
 \H_1 \vee \H_2& := \H_1 + \H_2 \, ,\label{eq:lattice2}
  \end{align}
shown on Fig.  \ref{fig:diagram}. Here  every embedding,  denoted by an arrow, is continuous and has dense range.
Taking conjugate duals,  it is easy to see that one has
\begin{align}
\H(R_G)^\times &= \H(R_G^{-1}) = \H + \H(G^{-1}),  \label{cup1}\\
\H(R_{G^{-1}})^\times &= \H(R_{G^{-1}}^{-1}) = \H + \H(G).   \label{cup2}
\end{align}
 In these relations, the r.h.s. is meant to carry the inductive norm (and topology) \cite[Sec.I.2.1]{pip-book}, so that both sides are in fact unitary equivalent, hence identified.

\medskip

At this stage, we return to the construction in terms of the closed unbounded operator $S$. We have to envisage two cases.
\medskip

\subsubsection*{(i)  An unbounded metric operator}
\label{subsec41}

We take as metric  operator $G_1= I+ S^\ast S$, which is unbounded, with $G_1 > 1$ and bounded inverse. Then the norm $\norm{G_1}{\cdot}$ is equivalent to the   norm  $\norm{R_{G_1}}{\cdot}$ on $D(G_1^{1/2}) = D(S)$, so that
$\H(G_1) = \H(R_{G_1})$ as vector spaces and thus also $\H(G_1^{-1}) = \H(R_{G_1}^{-1})$. On the other hand,  $G_1^{-1}$ is bounded.
Hence we get the triplet
\be\label{eq:tri>1}
\H(G_1) \; \subset\; \H \; \subset\; \H(G_1^{-1}) = \H({G_1})^\times.
\en
 Actually, the triplet \eqref{eq:tri>1}
 is the central part of the discrete scale of Hilbert spaces $V_{\G}$  built on the powers of  ${G_1^{1/2}}$.
This means that $V_{\G_1}:= \{\H_{n}, n \in \ZN \}$,
where $\H_{n} = D(G_1 ^{n/2}),  n\in \NN$, with a norm equivalent to the graph norm, and $ \H_{-n} =\H_{n}^\times$:
\be\label{eq:scale}
 \ldots\subset\; \H_{2}\; \subset\;\H_{1}\; \subset\; \H \; \subset\; \H_{-1} \subset\; \H_{-2} \subset\; \ldots
\en
Thus $\H_{1} =  \H({G_1}^{1/2} ) = D(S)$,  $\H_{2} =  \H({G_1} ) = D(S^\ast S)$, and $\H_{-1} =  \H({G_1}^{-1})$.
Then  $G_1^{1/2}$ is a unitary operator from $\H_1$ onto $\H$ and, more generally, from $\H_n$ onto $\H_{n-1}$.
In the same way,  $G_1$ is a unitary operator from $\H_n$ onto $\H_{n-2}$ and $G_1^{-1}$ s a unitary operator from $\H_n$ onto $\H_{n+2}$.

Moreover, one may add the end spaces of the scale, namely,
\be \label{eq:endscale}
\H_{\infty}({G_1} ):=\bigcap_{n\in \ZN} \H_n, \qquad \H_{-\infty}({G_1} ):=\bigcup_{n\in \ZN} \H_{n}.
\end{equation}
In this way, we get a genuine Rigged Hilbert Space:
\be
\H_{\infty}(G_1 ) \subset \H \subset \H_{-\infty}({G_1} ).
\en

In fact, one can go one more step. Namely, following \cite[Sec. 5.1.2]{pip-book}, we can use quadratic interpolation theory \cite{berghlof} and build a continuous scale of Hilbert spaces
$\H_{\alpha},   0\leq \alpha\leq 1$, between  $\H_{1}$  and $\H $, where $\H_{\alpha}=  D(G_1^{\alpha/2})$,  {with the graph norm  $\|\xi\|_{\alpha}^2 = \|\xi\|^2 + \|G_1^{\alpha/2}\xi\|^2$ or, equivalently, the norm
$\norm{}{(I+G_1)^{\alpha/2}\xi}^2$.
Indeed every $G_1^\alpha, \alpha\geq 0$, is   an unbounded metric operator.}

Next we define $\H_{-\alpha} =\H_{\alpha}^\times$ and iterate the construction to the full continuous scale $V_{\widetilde \G_1}:= \{\H_{\alpha}, \alpha \in \RN \}$.
Then, of course, one can replace $\ZN$ by $\RN$ in the definition \eqref{eq:endscale} of the end spaces of the scale.

\medskip

 In the general case,   $R_{G_1} = I+G_1 > 1$   is also an unbounded metric operator .
Thus we have
\be\label{eq:tri<1}
\H(R_{G_1}) \; \subset\; \H \; \subset\; \H(R_{G_1}^{-1}) = \H(R_{G_1})^\times,
\en
and we get another Hilbert-Gel'fand triplet.
Then one can repeat the construction  and obtain the Hilbert scale
 built on the powers of $R_{G_1}^{1/2}$.

 Now, if $S$ is injective, i.e. $N( S)= \{0\},$ then $|S|^2>0$ is also an unbounded metric operator.
Since $|S| >0, R_{S} = I+|S|^2 > 1$  and it is  another unbounded metric operator, with bounded inverse $R_{S}^{-1}$.
In both cases, one may build the corresponding Hilbert scale corresponding to the powers of $S$  or $R_S^{1/2}$.

At this stage, we have recovered the formalism based on metric operators that we have developed for the theory of quasi-Hermitian operators, in particular non-self-adjoint Hamiltonians, as encountered in the so-called $\P\T$-symmetric quantum mechanics. We refer to \cite{AT-metric1,AT-metric2,AT-metric3}  for a complete treatment.
However the case of an unbounded metric operator does not lead to many results, unless one considers  a \emph{quasi-Hermitian} operator \cite[Def.3.1]{AT-quHerm}.

\subsubsection*{(ii) A bounded metric operator}

We take as metric  operator $G_2= (I+ S^\ast S)^{-1}$, which is  bounded, with    unbounded inverse.

Since $G_2$ is \emph{bounded}, things simplify, because now $D(G_2) = \H$.
Similarly, one gets
$\H(R_{G_2^{-1}}) = \H(G_2^{-1})$ and $\H(R_{G_2^{-1}}^{-1}) = \H(G_2)$. So we are left with the triplet
\be
  \H(G_2^{-1}) \;\subset\; \H \;\subset\;  \H(G_2).
\label{eq:triplet}
\en
Then  $G_2^{1/2}$ is a unitary operator from $\H(G_2)$ onto $\H$ and from $\H$ onto $\H(G_2^{-1})$, whereas $G_2^{-1/2}$ is a unitary operator
$\H(G_2^{-1})$ onto $\H$ and from $\H$ onto $\H(G_2)$.

\medskip

{

\subsubsection*{(iii) A bounded metric operator with bounded inverse}

There is a third case, which is almost trivial. If the operator $S$ is bounded,   $G_1= I+ S^\ast S$ and $G_2= (I+ S^\ast S)^{-1}$ are both bounded metric operators, with bounded inverse. Then all nine \hs s in the lattice of Fig. \ref{fig:diagram} coincide    as vector spaces, with equivalent, but different, norms.

The advantage of this situation is that it leads to strong results on the similarity of two operators. As mentioned in \cite[Sec.3]{AT-PIPmetric}, up to unitary equivalence, one may always consider that the intertwining operator defining the similarity is in fact a metric operator. Let us briefly recall these notions.

Let $\H, \K$ be Hilbert spaces, $D(A)$ and $D(B)$ dense subspaces
of $\H$ and $\K$,  respectively,   $A:D(A) \to \H$, $B: D(B) \to \K$ two linear operators.
 A (metric) bounded operator $T:\H \to \K$ is called a   \emph{intertwining operator}  for $A$ and $B$ if
\begin{itemize}
\item[(i)] $T:D(A)\to D(B)$;
\item[(ii)] $BT\xi = TA\xi, \; \forall\, \xi \in D(A)$.
\end{itemize}
We say that $A$ and $B$ are \emph{similar}, and write $A\sim B$,
if there exists an intertwining operator $T$ for $A$ and $B$ with bounded inverse $T^{-1}:\K\to \H$, intertwining for $B$ and $A$.

A parallel definition (quasi-similarity) may be given in case the inverse $T^{-1}$ of the intertwining operator $T$ is not bounded.

From the relation $A\sim B$, there follows many interesting results about the respective spectra of $A$ and $B$, as described in detail in \cite{AT-PIPmetric}.
\medskip

\section{Transforming functions into frames by metric operators}

This section concerns the second main aim of our paper.
By Proposition \ref{prop 34}  and Theorem \ref{th_T^k}, we get the following result: if $\phi$ is a lower semi-frame of $\H$ and $D(C_\phi)$ is dense, then ${\sf T}_\phi^{-1/2}:\H\to \H$ and $\psi= {\sf T}_\phi^{-1/2} \phi$ is a Parseval frame for $\H_\phi = \H$. As we already remarked, ${\sf T}_\phi^{-1/2}:\H\to \H$ is a metric operator. We now want to relax the condition for $\phi$ of being a lower semi-frame. Thus we ask the following question:
for which weakly measurable functions $\phi:X\to \H$ does there exist a metric operator $G$ on $\H$ such that $G\phi$ is a frame? \footnote{Here and in the rest of paper,   writing $G\phi$  means implicitly that  $\phi_x\in D(G)$ for all $x\in X$.}
\smallskip

The motivation of this question is that in general one tries to pass from a less regular situation to a more regular, possibly in a smaller space.

In the next result we find some necessary or sufficient conditions for an answer to our question. As we are going to see, the recourse of the generalized frame operator is again useful.

\betheo
\label{th_funct_metric}
Let $\phi:X\to \H$ be a weakly measurable function with generalized frame operator ${\sf T}_\phi$. The following statements hold.
\begin{enumerate}
	\item[(i)] If $\phi$ is not total, then there does not exist a metric operator $G$ on $\H$ such that $G\phi$ is a frame and $G^{-1}$ is bounded.
	\item[(ii)] If $\phi$ is a Bessel mapping and not total, then there does not exist a metric operator $G$ on $\H$ such that $G\phi$ is a frame.
	\item[(iii)] If $D(C_\phi)$ is not dense, then there does not exist a metric operator $G$ on $\H$ such that $G\phi$ is a frame.
	\item[(iv)] If $D(C_\phi)$ is dense and $\phi$ is a lower semi-frame, then there exists a metric operator $G$ on $\H$ such that $G\phi$ is a Parseval frame. In particular, a possible choice is $G={\sf T}_\phi^{-1/2}$.
	\item[(v)] If $\phi$ is total, $D(C_\phi)$ is dense and $\phi_x\in R({\sf T}_\phi)$ for all $x\in X$, then there exists a metric operator $G$ on $\H$ such that $G\phi$ is a Parseval frame. In particular, a possible choice is $G={\sf T}_\phi^{-1/2}$.
\end{enumerate}
\entheo
\bdim
\begin{enumerate}
	\item[(i)] Suppose that there exists a metric operator $G$ on $\H$ such that $G\phi$ is a frame and $G^{-1}$ is bounded. By hypothesis, there exists $f\in \H, f\neq 0$ such that $\ip{\phi_x}{f}=0$ for a.e. $x\in X$. Since $G$ is self-adjoint and $G^{-1}$ is bounded, then $G^{-1}$ is defined on $\H$.
	Hence $\ip{G\phi_x}{G^{-1}f}=\ip{\phi_x}{GG^{-1}f}=0$  for a.e. $x\in X$, which implies that $G^{-1}f=0$, i.e. the contradiction $f=0$.

\item[(ii)] Suppose that there exists a metric operator $G$ such that $G\phi$ is a frame, and let ${\sf m}$ be a lower bound of $G\phi$. Then for all $f\in D(G)$
	\begin{equation}\label{G_Bessel}
	{\sf m}\|f\|^2\leq \int_X |\ip{f}{G\phi_x}|^2d\mu(x)=\int_X |\ip{Gf}{\phi_x}|^2d\mu(x)\leq {\sf M}\|Gf\|^2
	\end{equation}
	where ${\sf M}$ is an upper bound of $\phi$. Thus \eqref{G_Bessel} implies that $G^{-1}$ is bounded. By the previous point, we get then a contradiction.

\item[(iii)] Suppose that there exists a metric operator $G$ on $\H$ such that $G\phi$ is a frame. Then
	$
	\int_X |\ip{f}{G\phi_x}|^2 d\mu(x)<\infty$  for all $f\in \H.
	$
	In particular, for $f\in D(G)$ we have $$
	\int_X |\ip{Gf}{\phi_x}|^2 d\mu(x)<\infty,$$ so $R(G)\subseteq D(C_\phi)$. This means that $N(G)\supseteq D(C_\phi)^\perp\neq \{0\}$ which contradicts the property that $G$ is a metric operator.

	\item[(iv)] This follows by Proposition \ref{prop 34}.
	\item[(v)] The statement can be proved with a similar argument to that of Proposition \ref{prop 34} (Proposition \ref{prop32}), thus we give only a sketch of the proof.  First of all, we note that ${\sf T}_\phi$ has domain dense in $\H$, $D({\sf T}_\phi)\subseteq D({\sf T}_\phi^{1/2})$ and also that ${\sf T}_\phi$ is injective ($\phi$ is total). Let $f\in D({\sf T}_\phi)$, then
	\beano
	\int_{X} |\ip{f}{{\sf T}_\phi^{-1/2} \phi_x}|^2 \ud\mu(x)
	= \| {\sf T}_\phi^{1/2} {\sf T}_\phi^{-1/2} f \|^2= \|f \|^2.
	\enano
	Now, a standard argument of density concludes that $\psi$ is a Parseval frame of $\H$.
\end{enumerate}
\edim

\beex
Let $\{e_n\}_{n\in \NN}$ be an orthonormal basis, $\phi=\{e_1+e_n\}_{n\geq 2}$ and $\psi=\{e_n\}_{n\geq 2}$. Both $\phi$ and $\psi$ cannot be transformed into frames of $\H$ by a metric operator. Indeed, $D(C_\phi)=\{e_1\}^\perp$ and $\psi$ is a Bessel sequence but not total.
\enex

{
We will consider more examples (concerning, in particular, lower semi-frames) in the next section.  As an immediate consequence of Theorem \ref{th_funct_metric} we get the following result (compare with \cite[Corollary II.1]{kamuda}). We recall that two sequences $\{\phi_n\}_{n\in \NN},\{\psi_n\}_{n\in \mathbb{N}}$ are said {\em bi-orthogonal} if $\ip{\psi_n}{\phi_m}=\delta_{n,m}$, the Kronecker symbol.

\becor
Let $\phi:=\{\phi_n\}_{n\in \NN},\psi:=\{\psi_n\}_{n\in \mathbb{N}}$ be bi-orthogonal and total sequences of $\H$.
Then  $D(C_\phi)$ is dense. Moreover, if  ${\sf T}_\phi$ is the generalized frame operator of $\phi$, then $\phi_n\in R({\sf T}_\phi)\subset R({\sf T}_\phi^{1/2})$ for all $n\in \NN$ and $\{{\sf T}_\phi^{-1/2}\phi_n\}_{n\in \NN}$ is an orthonormal basis of $\H$.
\encor
\bdim
Then $D(C_\phi)$ is dense because it contains the total sequence $\psi$ and moreover $\psi_n\in D(T_\phi)$ and ${\sf T}_\phi \psi_n=\phi_n$ for all $n\in \mathbb{N}$ (which also gives $\phi_n\in R({\sf T}_\phi)\subset R({\sf T}_\phi^{1/2})$). Now by Theorem \ref{th_funct_metric}(iv) $\{{\sf T}_\phi^{-1/2}\phi_n\}_{n\in \NN}$ is a Parseval frame of $\H$, but since $\ip{{\sf T}_\phi^{-1/2}\phi_n}{{\sf T}_\phi^{-1/2}\phi_m}=\ip{\psi_n}{\phi_m}=\delta_{n,m}$  we conclude that $\{{\sf T}_\phi^{-1/2}\phi_n\}_{n\in \NN}$ is in particular an orthonormal basis of $\H$.
\edim

The problem about transforming functions in frames is still open. However, in the light of Theorem \ref{th_funct_metric} one may formulate a new version of the problem: given a weakly measurable function $\phi:X\to \H$, is it true that there exists  a metric operator $G$ on $\H$ such that $G\phi$ is a frame if and only if $\phi$ is total and $D(C_\phi)$ is dense?
}

\section{Examples}

In this final section, we exhibit several examples of lower semi-frames, mostly taken from our previous works.

\subsection*{(1) Sequences of exponential functions}

Let $\H=L^2(0,1)$ and $g\in L^2(0,1)$. A weighted exponential sequence of $g$ is $\mathcal{E}(g,b):=\{g_n\}_{n\in \ZN}=\{g(x)e^{2\pi i n b x}\}_{n\in \ZN}$ with $b>0$. In \cite[Remark 6.6]{corso2} it was proved that if $0<b\leq1$, then $\mathcal{E}(g,b)$ is a lower semi-frame if and only if $g$ is bounded away from zero, i.e. $|g(x)|\geq A$ for some $A>0$ and a.e. $x\in (0,1)$. Moreover,  by \cite[Corollary 6.5]{corso2}, the analysis operator of $\mathcal{E}(g,b)$ is densely defined and the generalized frame operator ${\sf T}_g$ of $\mathcal{E}(g,b)$ is the multiplication operator by $\frac{1}{b}|g|^2$, i.e.
$$
 D({\sf T}_g)=\{f\in L^2(0,1): f|g|^2\in L^2(0,1)\} \;\text{ and }\; {\sf T}_g f=\frac{1}{b}|g|^2 f, \;\;\forall f\in D({\sf T}_g).
$$
Thus, for $k\geq 0$ and $n\in \ZN$,
$({\sf T}_g^{-k}g_n)(x)=g(x)/|g(x)|^{-2k}e^{2\pi i n b x}, $
i.e.
${\sf T}_g^{-k}\mathcal{E}(g,b)=\mathcal{E}(g/|g|^{-2k},b)$. It is not difficult to check that the conclusions of Theorem \ref{th_T^k} hold.

\subsection*{(2) A Reproducing Kernel \hs}

Following \cite{at-reprodpairs2}, we consider   a reproducing kernel Hilbert space   of (nice) functions on
a measure space $(X, \mu)$, with kernel function $k_x, x\in X$, that is, $f(x)=\ip{f}{k_x}_K,\, \forall f\in\H_{K}$.
Choose the weight  function $m(x) >1$ and define the Hilbert scale $\H_k, \, k\in \ZN$, determined by the multiplication operator $Af(x) = m(x) f(x), \, \forall x\in X$. Thus we have,    for   $n \geq 1 $ $({\ov n} = -n),$
\be \label{quintet}
\H_{2n}  \subset \H_{n} \subset \H_K \subset  \H_{\ov n} \subset \H_{2\ov{n}} .
\en
Then define the measurable functions $\phi_x = k_x m^n(x), \psi_x = k_x m^{-n}(x)$, so that
 {$C_\psi : \H_K \to   \H_n,  \, C_\phi : \H_K \to  \H_{\ov n}$  continuously, and they are dual of each other. One has indeed
 $ \ip{\phi_x}{g}_K = \ov{g(x)}\,m^n(x) \in \H_{\ov n}$ and  $\ip{\psi_x}{g}_K = \ov{g(x)}\,m^{-n}(x) \in \H_{n}$, which implies duality.
 Thus $(\psi, \phi)$ is a reproducing pair with  $S_{\psi, \phi}=I$,  where $\psi$ is an upper semi-frame and $\phi$ a lower semi-frame.

Now we concentrate on the lower semi-frame $\phi$.
First we have  $D(C_\phi) = \H_{n}$, which is dense, so that $\H_\phi = \H$. 	Let us compute the sesquilinear form
\beano
	\Omega_\phi(f,g) &=& { \ip{ {C}_\phi f}{{C}_\phi g}}_K = \int_X f(x) m(x)^n\,\ov{g(x)}m(x)^n  \ud\mu(x),\;\; f,g\in D(C_\phi)
	\\
	&=&
	 \ip{{\sf T}_\phi^{1/ 2} f}{{\sf T}_\phi^{1/2} g}_K.
\enano
Therefore we have $ ({\sf T}^{1/ 2} _\phi f)(x) = f(x)\, m(x)^{n}$ and therefore $ ({\sf T}_\phi f)(x) = f(x)\, m(x)^{2n}$.
Hence $D({\sf T}_\phi) = \H_{2n} \subset  D(C_\phi) = \H_{n}$ and $D({\sf T}^{1/ 2}_\phi) = \H_{n} $, see \eqref{quintet} above.

Next ${\sf T}^{-1/ 2}_\phi $ is the operator of multiplication by $m^{-n}$, and indeed it is bounded in $\H$.
Finally, since $\phi$ is a lower semi-frame, $ \chi := {\sf T}^{-1/ 2}_\phi \phi$ is a frame in $\H$, by Proposition \ref{prop32}.

\subsection*{(3) Wavelets on the sphere}

The continuous wavelet   transform on the 2-sphere $\mathbb{S}^2$  has been analyzed in  \cite{wavsph}.
 For an axisymmetric (zonal) mother wavelet $\phi\in \H=L^2(\mathbb{S}^2,\ud\mu)$, define the family of spherical wavelets
  $$
  \phi_{\varrho,a}:=R_\varrho D_a\phi, \mbox{ where }(\varrho,a)\in X:= SO(3)\times \RN^+.
  $$
 Here, $D_a$ denotes the stereographic dilation operator  and $R_\rho$ the unitary rotation on $\mathbb{S}^2$.

Given the wavelet $\phi$,  it is known that the operator $S_\phi $ is diagonal in Fourier space (harmonic analysis on the 2-sphere reduces to expansions in spherical harmonics
 $Y_l^m, \, l \in \NN_0, m= -l, \ldots, l$), thus it is given by a Fourier multiplier   $\widehat{S_\phi f}(l,n)=s_\phi(l)\widehat f(l,n)$ with the symbol $s_\phi$ given by
  $$
  s_\phi(l):=\frac{8\pi^2}{2l+1}\sum_{|m|\leq l}\int_0^\infty\big|\widehat{D_a\phi}(l,m)\big|^2\frac{\ud a}{a^3},\ l\in\nN_0.
  $$
where $\widehat{D_a \phi}(l,m):= \ip {Y_l^m }  {D_a \phi} $ is the Fourier coefficient of $ D_a \phi$.

If one has ${\sf d}\leq s_\phi(l)  \leq {\sf c}$, for every $l \in \NN$, then the wavelet $\phi$ is admissible and a frame in $L^2(\mathbb{S}^2,\ud\mu)$.
However, it has been shown in \cite{wiaux} that the reconstruction formula converges under the weaker condition
 ${\sf d}\leq s_\phi(l)<\infty$ for all $l\in\nN_0$. In that case, $\phi$ is \emph{not}  admissible  and is a
 lower semi-frame, with  $S_\phi$ unbounded and densely defined. The domain of $S_\phi$ is the following:
$$
D(S_\phi =\{f  \in  L^2(\mathbb{S}^2,\ud\mu) : |\widehat f(l,n)|   \leq \frac{1} { s_\phi(l)}  |\widehat h(l,n)| , \text{  for some   } h\in L^2(\mathbb{S}^2,\ud\mu)\}.
$$
This domain contains, in particular, the set of band-limited functions, i.e. functions $f$ such that $\widehat f(l,n) = 0, \; \forall \,l \geq N_1$, for some $N_1 < \infty$, which is dense in $L^2$.
Since $S_\phi={D}_{\phi} {C}_{\phi}$, it follows that $D(S_\phi) \subset D(C_\phi) $, hence $ D(C_\phi)$ is dense as well and  $\H_\phi = \H$.

We proceed as in Case (2) and consider the sequilinear form
$$
	\Omega_\phi(f,g) = { \ip{ {C}_\phi f}{{C}_\phi g}}=
 \ip{{\sf T}_\phi^{1/ 2} f}{{\sf T}_\phi^{1/2} g}, \; f,g\in D(C_\phi).
$$
Since $ D(C_\phi)$ is dense, we get
${\sf T}_\phi:=|C_\phi|^2=C_\phi^* C_\phi$,
${\sf T}_\phi^{1/2} = |C_\phi|$ and $D({\sf T}_\phi) \subset D({\sf T}_\phi^{1/2}) = D(C_\phi).$

Finally, since $\phi$ is a lower semi-frame, $ \chi := {\sf T}^{-1/ 2}_\phi \phi$ is a frame in $\H$, by Proposition \ref{prop32}, as before.

\end{document}